\documentclass[11pt,twoside]{article}

\usepackage[english]{babel} \usepackage{amsmath} \usepackage{amssymb}
\usepackage{amsfonts} \usepackage[matrix,arrow]{xy}
\usepackage{pstricks} \usepackage[T1]{fontenc}
\usepackage{a4}\usepackage{theorem}

\newenvironment{rk}{\begin{quote} \normalfont\footnotesize {{\bf
        Remark} --} }{\end{quote}}

\newcommand{\pbullet}{$\bullet$~}

\newcommand{\sni}{\operatorname{sni}}
\newcommand{\Spec}{\operatorname{Spec}}
\newcommand{\Id}{\operatorname{Id}}
\newcommand{\ideng}[1]{\left\langle #1 \right\rangle}
\newcommand{\Perm}[1]{\mathcal{S}_{#1}} \newcommand{\QQ}{\mathbb Q}
\newcommand{\RR}{\mathbb R} \newcommand{\CC}{\mathbb C}
\newcommand{\PP}{\mathbb P} \newcommand{\OO}{\mathcal O}
\newcommand{\affecte}{\mathrel{:=}}
\newcommand{\scale}[2]{\genfrac{}{}{0pt}{}{#1}{#2}}
\newcommand{\proC}{{\mathbb P}\!\scale{1}{\CC}}
\newcommand{\proQ}{{\mathbb P}\!\scale{1}{\QQ}}
\newcommand{\per}{\operatorname{\mathcal S}}
\newcommand{\ant}{\operatorname{\mathcal A}}

\catcode`\à= \active \def à{\`a} \catcode`\À= \active \def À{\`A}
\catcode`\é= \active \def é{\'e} \catcode`\É= \active \def É{\'E}
\catcode`\è= \active \def è{\`e} \catcode`\È= \active \def È{\`E}
\catcode`\ù= \active \def ù{\`u} \catcode`\â= \active \def â{\^a}
\catcode`\Â= \active \def Â{\^A} \catcode`\ê= \active \def ê{\^e}
\catcode`\Ê= \active \def Ê{\^E} \catcode`\î= \active \def î{\^\i}
\catcode`\ô= \active \def ô{\^o} \catcode`\û= \active \def û{\^u}
\catcode`\ä= \active \def ä{\"a} \catcode`\ë= \active \def ë{\"e}
\catcode`\ï= \active \def ï{\"\i} \catcode`\ö= \active \def ö{\"o}
\catcode`\ü= \active \def ü{\"u} \catcode`\ç= \active \def ç{\c{c}}
\catcode`\Ç= \active \def Ç{\c{C}}

\catcode`\@=11 
\def \reserve#1/#2{\begingroup
     \dimen@\prevdepth \setbox0\vbox{#2}%
     \dimen@i=\dp0 \advance\dimen@i by\ht0 \ht0=0pt \dp0=\dimen@i
     \skip@\baselineskip \multiply\skip@ by#1 \advance\skip@
     by\baselineskip \global\dimen@i\wd0 \global\advance\dimen@i+5mm
     \nointerlineskip \nobreak
     \rightline{\vbox to0pt{\vskip\skip@\box0\vss}}%
     \nobreak \prevdepth\dimen@ \tolerance2000 \hfuzz2pt \hbadness2000
     \interlinepenalty\@M \def\par{\endgraf\nobreak}
     \everypar{\hangindent-\dimen@i \hangafter#1
       \everypar{\hangindent-\dimen@i \hangafter0 }}} \def
   \endreserve#1/{\ifnum#1=\z@ \endgroup \else
     \everypar{\hangindent-\dimen@i \hangafter-#1 \def\par{\endgraf
         \endgroup}} \fi}

\newcommand\scbf[1]{{%
    \usefont{T1}{cmr2}{bx}{sc} #1}}

\makeatletter \renewcommand{\@seccntformat}[1]{ \csname
  the#1\endcsname.\quad}
\renewcommand{\section}{\@startsection{section}{0}{0mm}{-\baselineskip}
  {0.5\baselineskip}{\normalfont\normalsize\centering\bfseries}}

\renewcommand{\subsection}{\@startsection{subsection}{1}{0mm}{\baselineskip}
  {-\fontdimen2\font plus -\fontdimen3\font minus -\fontdimen4\font}
  {\normalfont\normalsize\bfseries}}

\renewcommand{\subsubsection}{\@startsection{subsubsection}{1}{0mm}{\baselineskip}
  {-\fontdimen2\font plus -\fontdimen3\font minus -\fontdimen4\font}
  {\normalfont\normalsize\slshape}} \makeatother

\theoremstyle{plain} \theoremheaderfont{\bfseries\slshape}
\theorembodyfont{\slshape} \newtheorem{theo}{Theorem}[section]

\newtheorem{lem}[theo]{Lemma} \newtheorem{prop}[theo]{Proposition}
\newtheorem{fact}[theo]{Fact} 

\title{
  {\scbf{Computation of some Moduli Spaces of covers} \\
    {\small \bf{and explicit $\per_n$ and $\ant_n$ regular
        $\QQ(T)$-extensions with totally real fibers}}}}

\author{Emmanuel {\sc Hallouin} and Emmanuel {\sc
    Riboulet-Deyris}\thanks{Groupe de Recherche en Informatique et
    Mathématique du Mirail (G.R.I.M.M.), University of Toulouse~II,
    France.}}

\begin{document}
\maketitle

\begin{quote}
  \footnotesize {\bfseries Summary. ---} We study and compute an
  infinite family of Hurwitz spaces parameterizing covers of~$\proC$
  branched at four points and deduce explicit regular~$\mathcal{S}_n$
  and~$\mathcal{A}_n$-extensions over~$\QQ(T)$ with totally real
  fibers.
\end{quote}

\section*{Introduction}

In this paper, we study a family of covers of the projective line
suggested to us by Gunter Malle, namely those covers of even degree~$n
\geq 6$, ramified over four points, with monodromy $\per_n$ and having
branch cycle description~$\mathbf{C} = (C_1, C_2, C_3, C_4)$ of type:
$$
\left( (n-2),3,2^{\frac{n-2}{2}},2^{\frac{n}{2}} \right)
$$
Malle suspected the Hurwitz curves have genus zero for every~$n
\geq 6$ and some covers in the family have totally real fibers. A
similar family was suggested and extensively studied by Dèbes and
Fried in~\cite{Debes_Fried_Non_Rigid}. Unfortunately, their Hurwitz
spaces happen to have a quadratic genus in~$n$ and only provide the
expected regular extensions for degrees~$5$ and~$7$
(see~\cite{Debes_Fried_Non_Rigid}, Theorem~$4.11$).  Their work uses
braid action formulae (see~\cite{Fried_Volklein}) and complex
conjugation action formulae (see~\cite{Fried_Debes_Real_Res_Field},
Proposition~$2.3$).

In this paper, we first follow the lines of Dèbes and Fried method and
show that Malle's expectations were right. Then the second half of the
paper is devoted to the explicit calculation of the concerned
universal family of covers. To this end, we use an explicit version of
Harbater's deformation techniques (\cite{Harbater_deformation}) as
proposed in~\cite{JMC_Granboulan, JMC_Tools}. As far as we know, it is
the first time such an advanced method is used for computing an
infinite family of covers. We present numerical results, obtained with
{\sc Magma}, showing the efficiency of the proposed method compared to
the classical ones involving Buchberger algorithm. Indeed our
computation reduces to solving {\slshape linear systems}.

In the first section, we recall results from the theory of Hurwitz
spaces and (non)-rigidity methods developed in~\cite{Fried_Volklein,
  Volklein, Fulton_Hur,Malle_Matzat}.  The second section is devoted
to the combinatoric study of our family and arithmetical consequences
of it using the method of Dèbes and Fried. In the third section, we
show the existence of totally real $\Perm{n}$ and~$\mathcal{A}_n$
residual $\QQ$-extensions in that family. Related to this totally real
specialization, we find a very special point of the boundary of our
Hurwitz space that shows very useful when computing an explicit model.
This is done in the last section, by a deformation method.

We notice that a fallout of our construction is the existence of
totally real $\QQ$-extensions with Galois group~$\mathcal{S}_n$
and~$\mathcal{A}_n$. However a simpler and astute construction can be
found in~\cite{Mestre}.

We thank Jean-Marc Couveignes for numerous extremely helpful
discussions about this work.

\section{General framework}

Associated to our family of covers, there is a coarse moduli space
called Hurwitz space, denoted~$\mathcal{H}'_4(\Perm{n}, \mathbf{C})$.
It is a quasi-projective regular (not a priori connected) variety
over~$\overline{\QQ}$ with the following properties
(see~\cite{Fried_Volklein, Volklein, Fulton_Hur}):

\smallskip

\pbullet Since any conjugacy class of~$\Perm{n}$ is
rational,~$\mathcal{H}'_4(\Perm{n}, \mathbf{C})$ is defined over~$\QQ$.

\pbullet Let~$F_4$ be the configuration space of $4$ points, e.g.
$(\proC)^4 \setminus \mathcal{D}$ where $\mathcal D$ denotes the
discriminant variety. The map:
$$
\xymatrix@R=7pt@C=4pt{
  *[l]{\phi :} & \mathcal{H}'_4(\Perm{n}, \mathbf{C}) \ar[r] & F_4 (\proC)\\
  & h \quad \ar@{|->}[r] & (z_1,z_2,z_3,z_4)\\
  }
$$
where~$z_1,z_2,z_3,z_4$ are the branched points (in the given
order) of the cover corresponding to the point~$h$, is a finite étale
morphism defined over~$\QQ$.

\pbullet Since~$\Perm{n}$ has no center, the covers in our family have
no automorphism so the moduli space~$\mathcal{H}'_4(\Perm{n},
\mathbf{C})$ is a fine one and for any~$h \in \mathcal{H}'_4(\Perm{n},
\mathbf{C})$ the associated cover~$p_h$ can be defined over~$\QQ(h)$
the field of definition of the point~$h$.

\smallskip

As in \S$4.2$ of~\cite{Debes_Fried_Non_Rigid}, rather than looking at
the full moduli space, we concentrate on a curve in it.  Let us fix
three points~$z_1, z_2, z_3 \in \PP^1 (\QQ)$ and consider the
curve~$\mathcal{H}'_{(z_1, z_2, z_3)}$ obtained by the pullback:
$$
\xymatrix{ \mathcal{H}'_{(z_1, z_2, z_3)} \ar[r] \ar[d]_{\varphi} &
  \mathcal{H}'_4(\Perm{n}, \mathbf{C}) \ar@{->}[d]^{\phi} \\
  \PP^1_{\CC} \setminus \{z_1, z_2, z_3\} \ar[r]^{\qquad i} & F_4
  (\proC) }
$$
(the lower horizontal map~$i$ is~$z \mapsto (z_1, z_2, z_3, z)$).
If the three fixed points are rational, all the maps are defined
over~$\QQ$ and the curve~$\mathcal{H}'_{(z_1, z_2, z_3)}$ is also
defined over~$\QQ$.

\section{Combinatoric study of the Hurwitz curve}

In this section, we gather the more information we can about the cover
$\varphi : {\mathcal H}'_{(z_1, z_2,z_3)} \longrightarrow
{\proC}\smallsetminus \{ z_1 , z_2 , z_3 \}$ from its combinatoric
description, namely the monodromy, the
ramification type, the connectivity and the genus. This will be deduced,
from an explicit enumeration of the Nielsen class associated
to~$\mathbf{C}$ (section~\ref{combinatoiredelaclassedenielsen}) and
from the action of braids on it (section~\ref{s_braiding_action}).

\subsection{Nielsen classes description. \normalfont ---} \label{combinatoiredelaclassedenielsen}

Let us fix the branched points~$\underline{z} = (z_1,z_2,z_3,z_4) \in
F_4 (\proC)$ and an homotopic base of $\proC \setminus \{
z_1,z_2,z_3,z_4\}$. Using the topological classification of covers,
elements of the fiber~$\phi^{-1}(\underline{z})$ are in bijection
with~$\sni^{\rm ab}(\Perm{n}, \mathbf{C})$ the strict absolute Nielsen
class of type~$(\Perm{n}, \mathbf{C})$, that is:
$$
\displaystyle{ \sni^{\rm ab}(\Perm{n}, \mathbf{C}) = \left\{
    (\sigma_1, \sigma_2, \sigma_3, \sigma_4) \in (\Perm{n})^4, \;
    \left\{\begin{array}{l}
        \sigma_1 \sigma_2 \sigma_3 \sigma_4 = 1, \; \sigma_ i \in C_i \; \forall i \\
        \ideng{\sigma_1, \sigma_2, \sigma_3, \sigma_4} = \Perm{n}\\
                 \end{array}
               \right.  \right\} \Biggr/ \sim }
           $$
           where~$(\sigma_1, \sigma_2, \sigma_3, \sigma_4) \sim
           (\sigma'_1, \sigma'_2, \sigma'_3, \sigma'_4)$ means that
           there exists~$\tau \in \Perm{n}$ such that~$\sigma'_i =
           \tau \sigma_i \tau^{-1}$ for all~$1 \leq i \leq 4$. We
           first enumerate the Nielsen class.
           
           Since any two $(n-2)$-cycles are $\Perm{n}$-conjugate,
           every element of~$\sni^{\rm ab}(\Perm{n}, \mathbf{C})$ has
           a representative with~$\sigma_1 = (1, \ldots, n-2)$.
           Conjugating by a power of~$\sigma_1$, we also assume
           that~$\sigma_2 = (n-2, k, l)$. We now distinguish three
           cases.
           
           \pbullet {\sl The case~$\{k,l\} = \{n-1, n\}$.} Conjugating
           by a power of~$\sigma_1$ and by~$(n-1,n)$, every element in
           that class has a unique representative with~$\sigma_1 = (1,
           \ldots, n-2)$, and~$\sigma_2 = (n-2,n-1,n)$. So~$\sigma_1
           \sigma_2 = (1, \ldots, n)$ and the enumeration reduces to
           finding all the permutations~$(\sigma_4, \sigma_3) \in C_4
           \times C_3$ such that~$\sigma_4 \sigma_3 = (1, \ldots, n)$
           and~$\ideng{\sigma_1, \sigma_2, \sigma_3, \sigma_4} =
           \Perm{n}$.  We need a lemma used several times further. It
           deals with relations in the diedral group~$D_m$.

\begin{lem} \label{tuilage}
  Let~$m \leq n$ be even and let~$c$ be an $m$-cycle of~$\Perm{n}$.
  There is a bijection between non trivial cycles of~$c^{\frac{m}{2}}$
  (i.e. sets of the form~$\{x, c^{\frac{m}{2}}(x)\}$ where~$x$ belongs
  to the support of~$c$), and the decompositions~$c = \sigma \tau$
  with~$\sigma$ a product of~$\frac{m}{2}$ transpositions and~$\tau$ a
  product of~$\frac{m}{2} - 1$ transpositions. Therefore there are
  exactly~$\frac{m}{2}$ such decompositions of~$c$.
\end{lem}

\noindent
{\slshape Proof}. --- Let~$x$ be an element of the support of~$c$ then
one can verify that:
\begin{equation} \label{sigma_tau}
c =
\underbrace{
\prod_{i = 1}^{\frac{m}{2}} \left( c^{1-i}(x), c^i(x) \right)
}_{\sigma_{c, x}}
\underbrace{
\prod_{j = 1}^{\frac{m}{2} - 1} \left( c^j(x), c^{-j}(x) \right)
}_{\tau_{c, x}}
\end{equation}
This is a decomposition associated to the set~$\{x,
c^{\frac{m}{2}}(x)\}$. Two such decompositions associated to~$x$
and~$y$ are equal if and only if~$\{x, c^{\frac{m}{2}}(x)\} = \{y,
c^{\frac{m}{2}}(y)\}$. On the other hand, let~$c = \sigma \tau$ be a
decomposition as in the lemma then one can show that it can be written
as in~(\ref{sigma_tau}) by considering an element~$x$ of the support
of~$\sigma$ not belonging to the support of~$\tau$.\hfill $\square$

The group~$D_m$ is known to be the group of isometries of a regular
polygon with $m$~vertices.  We can explain the previous relation
geometrically: for each vertex~$x$, the rotation (i.e.~$c$) inducing
an $m$-cycle on the vertices equals the composition of the unique
reflection permuting the vertex~$x$ and its successor
(i.e.~$\sigma_{c,x}$) with the unique reflection fixing~$x$
(i.e.~$\tau_{c,x}$).

Back to the enumeration of the Nielsen classes in the case
where~$\{k,l\} = \{n-1,n\}$, the lemma~\ref{tuilage} shows that there
are exactly~$\frac{n}{2}$ such elements; this subset is denoted
Class~$\mathcal{A}$ in table~\ref{Tableau_recapitlatif}.

\pbullet {\em The case~$\#(\{k,l\} \cap \{n-1, n\}) = 1$.}
Conjugating by a power of~$\sigma_1$ and by~$(n-1, n)$, every class
has a unique representative with~$k \in \{1, \ldots, n-3\}$,~$\sigma_1
= (1, \ldots, n-2)$ and~$\sigma_2 = (k,n-2,n-1)$.  So~$\sigma_1
\sigma_2 = (1, \ldots, k)(k+1,\ldots, n-1)(n)$; but~$\sigma_1 \sigma_2
= \sigma_4 \sigma_3$ has necessarily an even number of fixed points
(because if~$x$ is fixed by~$\sigma_4 \sigma_3$, so is~$\sigma_4(x)
\not= x$) therefore we must have~$k = 1$ and~$\sigma_1 \sigma_2 = (2,
\ldots, n-1)$. We stress that if~$x$ is fixed by~$\sigma_4 \sigma_3$
then~$\sigma_3(x) = \sigma_4(x) \not= x$ and~$\sigma_4(x)$ is also
fixed by~$\sigma_4 \sigma_3$. Therefore~$\sigma_3$ and~$\sigma_4$
share the transposition~$(1,n)$ and the rests of the decompositions
are given by lemma~\ref{tuilage} with~$c = (2, \ldots, n-1)$.  These
elements form the Class~$\mathcal{B}$ in
table~\ref{Tableau_recapitlatif}.

\pbullet {\em The case~$\{k,l\} \subset \{1, \ldots, n-3\}$.} In that
case~$\sigma_1 = (1, \ldots, n-2)$ and~$\sigma_2 = (k,l,n-2)$, then:
$$
\begin{cases}
  \sigma_1 \sigma_2 = (1, \ldots, k,l+1, \ldots, n-2,k+1, \ldots, l)(n-1)(n) & \text{if~$k < l$} \\
  \sigma_1 \sigma_2 = (1, \ldots, l)(l+1, \ldots, k)(k+1, \ldots,
  n-2)(n-1)(n) & \text{if~$k > l$}
\end{cases}
$$
So~$\sigma_1 \sigma_2 = \sigma_4 \sigma_3$ fixes~$n-1$ and~$n$.
Note that~$\sigma_3$ and~$\sigma_4$ cannot share the
transposition~$(n-1, n)$, because if so~$\{n-1, n\}$ would be
an~$\ideng{\sigma_1, \sigma_2, \sigma_3, \sigma_4}$-orbit; this would
contradict the transitivity.  Therefore~$\sigma_1 \sigma_2 = \sigma_4
\sigma_3$ has at least four fixed points:~$n-1, \sigma_4(n-1), n$
and~$\sigma_4(n)$. Then we must have~$k > l$ and, two of the
lengths~$l$,~$k-l$ and~$n-2-k$ equal~$1$. There are three
possibilities~$\sigma_2 = (n-3, 1, n-2)$ or~$\sigma_2 = (2, 1,n-2)$
or~$\sigma_2 = (n-3, n-4, n-2)$.  which are all together conjugate
under a power of~$\sigma_1$. In conclusion every element in this class
has a unique representative such that~$\sigma_1 = (1, \ldots, n-2)$
and~$\sigma_2 = (n-2,n-3,n-1)$.  Then~$\sigma_1 \sigma_2 = (1, \ldots,
n-4)$ and in the support of~$\sigma_3$ and~$\sigma_4$, we
find~$(n,n-2)(n-1,n-3)$ or~$(n,n-3)(n-1,n-2)$ which are conjugate
under~$(n-1,n)$. Again, lemma~\ref{tuilage} for~$c = (1, \ldots, n-4)$
concludes.  These elements form the Class~$\mathcal{C}$ in
table~\ref{Tableau_recapitlatif}.

We note that for every $4$-tuple, the group generated by the four
permutations is~$\mathcal{S}_n$.

The whole enumeration can be find in table~\ref{Tableau_recapitlatif}.
Note that the three pointed classes have the following cardinalities:
$$
\# \mathcal{A} = \frac{n}{2}, \qquad \# \mathcal{B} = \frac{n}{2} -
1 \qquad \text{and} \qquad \# \mathcal{C} = \frac{n}{2} - 2
$$
therefore:
$$
\# \sni^{\rm ab}(\Perm{n}, \mathbf{C}) = 3\left( \frac{n}{2} -
  1\right)
$$
Concerning the Hurwitz cover, we have shown that:
\begin{fact} \label{degree_Hurwitz_cover}
  The degree of the Hurwitz cover~$\phi$ (or~$\varphi$)
  equals~$3\left( \frac{n}{2} - 1\right)$.
\end{fact}

\begin{table}[!ht]
$$
\begin{array}{|c|l|}
\hline
\text{Class~$\mathcal{A}$} &
a_i = \left[
        (1, \ldots, n-2), (n-2,n-1,n), \tau_{(1, \ldots, n), i}, \sigma_{(1, \ldots, n), i}
      \right] \\
 & \text{with~$1 \leq i \leq \frac{n}{2}$}\\
\hline
\text{Class~$\mathcal{B}$} &
b_i = \left[
        (1, \ldots, n-2), (1,n-2,n-1), \nu\tau_{(2, \ldots, n-1), i},
                                       \nu\sigma_{(2, \ldots, n-1), i} \right] \\
 & \text{with~$\nu = (1,n)$ and~$2 \leq i \leq \frac{n}{2}$} \\
\hline
\text{Class~$\mathcal{C}$} &
c_i = \left[
        (1, \ldots, n-2), (n-2,n-3,n-4), \nu\tau_{(1, \ldots, n-4), i},
                                         \nu\sigma_{(1, \ldots, n-4), i} \right] \\
 & \text{with~$\nu = (n,n-2)(n-1,n-3)$ and~$1 \leq i \leq \frac{n}{2} - 2$} \\
\hline
\end{array}
$$
\caption{\label{Tableau_recapitlatif} The Nielsen classes
(same notations as in lemma~\ref{tuilage})}
\end{table}

\subsection{Braiding action. \normalfont ---} \label{s_braiding_action}

In this paragraph, we compute the action of braids on the Nielsen
class given in table~\ref{Tableau_recapitlatif}.

\subsubsection{Generator of the braid group and braiding action. \normalfont ---} Let
us denote by ${\mathcal H}_4 (\Perm{n}, \mathbf{C})$ the Hurwitz
space parameterizing the same set of isomorphism classes of covers as
$\mathcal{H}'_4 (\Perm{n}, \mathbf{C})$ but without ordering the branch
points. This space can be endowed with a topology which is constructed
in the same way as the one of $\mathcal{H}'_4 (\Perm{n}, \mathbf{C})$
(see \cite{Fried_Volklein} or \cite{Fulton_Hur}). On the other hand,
it maps onto~$C_4 (\proC)$, the quotient of~$F_4 (\proC)$ by the
action of~$\per_4$ on the coordinates:
\begin{equation}
\vcenter{
\xymatrix@R=14pt@C=12pt{
\mathcal{H}'_4 (\Perm{n}, \mathbf{C}) \ar[r] \ar[d]^{\phi} & {\mathcal H}_4 (\Perm{n}, \mathbf{C}) \ar[d]^{\phi'} \\
F_4 (\proC) \ar[r] & C_4 (\proC) }}
\end{equation}
The fundamental group of~$C_4 (\proC)$ is the Hurwitz braid group of
index $4$. It possesses a classical presentation (see \cite{hansen} or
\cite{birman}):
$$
\left \langle Q_1 , Q_2 , Q_3 \;
\begin{array}{|l}
Q_1 Q_3 = Q_3 Q _1\\ 
Q_1 Q_2 Q_1 = Q_2 Q_1 Q_2 \text{ and }  Q_2 Q_3 Q_2 = Q_3 Q_2 Q_3 \\
Q_1 Q_2 Q^{2}_3 Q_2 Q_1 =1 \: \text{\small (sphere's relation)}\\
\end{array} \right \rangle
$$

Denoting the fiber of $\phi'$ by ${\operatorname{ni}}^{\rm
  ab}(\Perm{n}, \mathbf{C})$ the ``unordered'' Nielsen class
associated to inertia's $4$-tuple $\mathbf{C}$, i.e. the quotient of
$\sni^{\rm ab}(\Perm{n}, \mathbf{C})$ by the action of $\per_4$ on the
coordinates, we have:
\begin{prop}[Braiding action formula]
  For all $i = 1,2,3$, and for all $[\sigma_1 , \sigma_2 , \sigma_3 ,
  \sigma_4]$ in ${\operatorname{ni}}^{\rm ab}(\Perm{n}, \mathbf{C})$,
  the monodromy (right) action of the braid~$Q_i$ for the previous
  generating system is given by the formula:
  $$[\sigma_1 , \sigma_2 , \sigma_3 , \sigma_4].Q_i = [\ldots ,
  \sigma_i \sigma_{i+1}\sigma^{-1} _i , \sigma_i, \ldots]$$
  where
  $\sigma_i \sigma_{i+1}\sigma^{-1} _i$ is the $i$-th coordinate.
\end{prop}
Generally speaking the diagram ($2$) permits us to express the
monodromy of $\phi$ according to the one of $\phi'$. This can be done
explicitly because a presentation of the fundamental group of $F_4
(\proC)$ can be expressed in term of the one of $C_4 (\proC)$. We just
give here the generators and refer to \cite{hansen} for a complete
system of relations:
$$
t_{1,2} = {Q_1}^2, \qquad t_{2,3} = {Q_2}^2, \qquad t_{1,3} =
{Q_1}{Q_2}^2{Q_1}^{-1}
$$

\subsubsection{Monodromy action for the cover $\varphi$. \normalfont ---} Let us
recall that we just defined~$\varphi$ in the section $1$ to be the
pullback of $\phi$ by the monomorphism denoted by $ i: \proC
\smallsetminus \{ z_1 , z_2 , z_3 \} \hookrightarrow {\mathcal U}^4$.
We now want to study the monodromy of this cover.

Let us choose $z_1 < z_2 < z_3$ on the real line $\PP^1(\RR)$ and let
us define the following ``homotopic base'' of $\proC \smallsetminus
\{z_1 , z_2 , z_3\}$:
\begin{center}
\begin{pspicture}(-5,-0.75)(5,1)
  \psline[linewidth=0.01, linestyle=dashed](-5,0)(5,0)
  \psdots(-3,0)(-1,0)(1,0)(3,0) \rput(-3.1,-0.2){$z_1$}
  \rput(-1.1,-0.2){$b$} \rput(0.9,-0.2){$z_2$} \rput(2.9,-0.2){$z_3$}
  \rput(5,-0.4){$\PP^1(\RR)$}
  \psarc[linewidth=0.03,arrowsize=8pt
  2,arrowlength=0.6,arrowinset=0.75]{->}(-3,0){0.5}{300}{300}
  \psline[linewidth=0.03](-1,0)(-2.5,0) \rput(-1.75,-0.2){$\gamma_1$}
  \psarc[linewidth=0.03,arrowsize=8pt
  2,arrowlength=0.6,arrowinset=0.75]{->}(1,0){0.5}{300}{300}
  \psline[linewidth=0.03](-1,0)(0.5,0) \rput(-0.25,-0.2){$\gamma_2$}
  \pscurve[linewidth=0.03](-1,0)(1,0.75)(2.5,0)
  \psarc[linewidth=0.03,arrowsize=8pt
  2,arrowlength=0.6,arrowinset=0.75]{->}(3,0){0.5}{300}{300}
  \rput(-0.1,0.75){$\gamma_3$}
\end{pspicture}
\end{center}
Then, as explained in details in Theorem~$4.5$
of~\cite{Debes_Fried_Non_Rigid}, for an adapted generating system of
braids\footnote{In our situation the $3$ points $z_1$, $z_2$ and $z_3$
  are on the real line, so the choice of $Q_1$, $Q_2$, $Q_3$ is the
  standard one (see \cite{hansen} for example). In general the $3$
  generators $Q_i$ must be precisely given. We point out that this
  choice just depends on a given path through~$z_1$,~$z_2$,~$z_3$ in this
  order.}  $Q_1$, $Q_2$, $Q_3$ on $C_4 (\proC)$, we can compute the
morphism $i_{\ast}$ induced by $i$ on the respective homotopic groups
in term of those two sets of generators:
$$
i_{\ast} (\gamma_1) = Q^{2} _1= t_{1,2}, \quad i_{\ast} (\gamma_2)
= Q^{2} _2 = t_{2,3}, \quad i_{\ast} (\gamma_3) = Q^{-1} _2 Q^{2} _3
Q_2 = (t_{1,2} t_{2,3})^{-1}
$$
With these formulas, the computation of the monodromy of $\varphi$
is deduced from the one of $\phi$. We summarize in the next:
\begin{prop}
  Using the notations of the table~$1$, the monodromy action for the
  cover $\varphi$ is:
  
  $\bullet$ for the path $\gamma_1$:
  $$(a_{\frac{n}{2}} , a_{\frac{n}{2}-1} , \ldots ,
  a_1)(b_{\frac{n}{2}}, b_{\frac{n}{2}-1}, \ldots ,
  b_2)(c_{\frac{n}{2}-2}, c_{\frac{n}{2}-3} , \ldots , c_1)$$
  
  $\bullet$ for the path $\gamma_2$:
\begin{equation*}
\begin{split}
  & (a_{\frac{n}{2} - 2} , b_\frac{n}{2}, b_{\frac{n}{2}-1},
  a_\frac{n}{2}, c_{\frac{n}{2}-2})(a_1 , c_{\frac{n}{2}-3})(a_2 ,
  c_{\frac{n}{2}-4})\ldots (a_{\frac{n}{2}-3} , c_{1})\\
  &(a_{\frac{n}{2}-1})
\begin{cases}
  (b_2, b_{\frac{n}{2}-2})(b_3, b_{\frac{n}{2}-3})\ldots(b_{\frac{n}{4}-1}, b_{\frac{n}{4}+1})(b_{\frac{n}{4}}) & \text{ if $4 \mid n$,}\\
  (b_2, b_{\frac{n}{2}-2})(b_3, b_{\frac{n}{2}-3})\ldots(b_{\frac{n-2}{4}}, b_{\frac{n+2}{4}}) & \text{ if $4 \nmid n$,}\\
\end{cases}
\end{split}
\end{equation*}

$\bullet$ for the path $\gamma_1 . \gamma_2$:
\begin{equation*}
\begin{split}
  (& a_{\frac{n}{2}-1}  , a_{\frac{n}{2}-2} , b_{\frac{n}{2}-1})(c_1 , a_{\frac{n}{2}-4})(c_2 , a_{(\frac{n}{2}-5)}) \ldots (c_{\frac{n}{2}-4} , a_{1})(c_{\frac{n}{2}-3} , a_\frac{n}{2})(c_{\frac{n}{2}-2} , a_{\frac{n}{2} - 3})\\
  &
\begin{cases}
  (b_{\frac{n}{2}} , b_{\frac{n}{2}-2})(b_{\frac{n}{2}-3} , b_2)(b_{\frac{n}{2}-4} , b_3)\ldots (b_{\frac{n}{4}+1} , b_{\frac{n}{4} - 2})(b_{\frac{n}{4}} , b_{\frac{n}{4}-1}) & \text{if $4 \mid n$,}\\
  (b_{\frac{n}{2}} , b_{\frac{n}{2}-2})(b_{\frac{n}{2}-3} , b_2)\ldots
  (b_{\frac{n+2}{4}} , b_{\frac{n-2}{4}-1})(b_\frac{n-2}{4})& \text{if
    $4 \nmid n$,}
\end{cases}     
\end{split}
\end{equation*}
\end{prop}

\subsection{Ramification in the Hurwitz curve. \normalfont ---}\label{Hurwitz_curve_ramification}

In conclusion, the Hurwitz curve~$\mathcal{H}'$ is a cover
of~$\PP^1_{\CC}$ of degree~$3 \left(\frac{n}{2} - 1\right)$ ramified
over three points, say~$z_1, z_2, z_3 \in \PP^1 (\QQ)$ with
ramification type described in
table~\ref{Tab_Hurwitz_curve_ramification}.

\begin{table}[!ht]
$$
\begin{array}{|c|}
\hline
\text{if~$4 \mid n$} \\
\hline
\xymatrix@C=0pt@R=8pt{
\frac{n}{2} & \frac{n}{2} - 1 & \frac{n}{2} - 2 & 5 & 1^2 & 2^{\frac{3n}{4}-5} & 3 &           & 2^{\frac{3n}{4} - 3} \\
            & z_1 \ar@{-}[ul]\ar@{-}[u]\ar@{-}[ur]      &  &               & z_2 \ar@{-}[ul]\ar@{=}[u]\ar@{=}[ur]      &   &   & z_3 \ar@{-}[ul]\ar@{=}[ur]&
} \\
\hline
\text{if~$4 \nmid n$} \\
\hline
\xymatrix@C=0pt@R=8pt{
\frac{n}{2} & \frac{n}{2} - 1 & \frac{n}{2} - 2 & 5 & 1 & 2^{\frac{3(n-6)}{4}} & 3 & 1 & 2^{\frac{3n-14}{4}} \\
& z_1 \ar@{-}[ul]\ar@{-}[u]\ar@{-}[ur] & & & z_2 \ar@{-}[ul]\ar@{-}[u]\ar@{=}[ur] & & & z_3 \ar@{-}[ul]\ar@{-}[u]\ar@{=}[ur]&
²} \\
\hline
\end{array}
$$
\caption{\label{Tab_Hurwitz_curve_ramification} Ramification in the Hurwitz curve}
\end{table}

\begin{fact}\label{defined_over_Q}
  The Hurwitz curve ${\mathcal H}'_{(z_1,z_2,z_3)}$ satisfies:
  \begin{itemize}
  \item it is irreducible;
  \item it is of genus zero and $\QQ$-isomorphic to $\proQ$.
  \end{itemize}
  So, our family contains infinitely many covers defined over~$\QQ$.
\end{fact}

\noindent
{\slshape Proof}. --- The irreducibility comes from the transitivity
of the braiding action. The Riemann-Hurwitz formula shows that the
genus of~$\mathcal{H}'_{(z_1,z_2,z_3)}$ is zero. We can also note that,
for example, the point of ramification index~$5$ must be a rational
one (this ramification index is isolated).  Being defined over~$\QQ$,
of genus zero and with a rational point,~$\mathcal{H}'_{(z_1,z_2,z_3)}$
is necessarily $\QQ$-isomorphic to~$\PP^1_\QQ$. In particular, there
are covers in our family defined over $\QQ$.\hfill $\square$

\section{Existence of totally real~$\Perm{n}$ and~$\mathcal{A}_n$-extensions}

There is still a question left: does our family contain elements
defined over~$\QQ$ with totally real fibers?

In this section, we use complex conjugation action on fibers as
describe by Dèbes and Fried (see theorem~$2.4$
of~\cite{Fried_Debes_Real_Res_Field} or proposition~$2.3$
of~\cite{Debes_Fried_Non_Rigid}) and prove adapted formulae to our
family and to our choice of homotopic basis.

\subsection{Covers in the family with totally real fiber. \normalfont ---} \label{s_totally_real_fibers}

We consider a finite algebraic cover~$p : \mathcal{C} \rightarrow
\PP^1_\CC$ ramified over four ordered real points~$z_2 < z_3 < z_4 <
z_1 \in \PP^1(\RR)$, we fix~$z_0 \in \PP^1(\RR)$ a real base point
between~$z_3$ and~$z_4$ and we denote by~$F$ the fiber~$p^{-1}(z_0)$.

The complex conjugation will play a crucial role; let denote this
conjugation by a bar; for example~$\bar{p} : \bar{\mathcal{C}}
\rightarrow \PP^1_\CC$ is the cover obtained from the first one by
complex conjugation and~$\bar{F}$ is its fiber above~$z_0$.  Let~$c :
F \rightarrow \bar{F}$ be the bijection induced by the complex
conjugation.

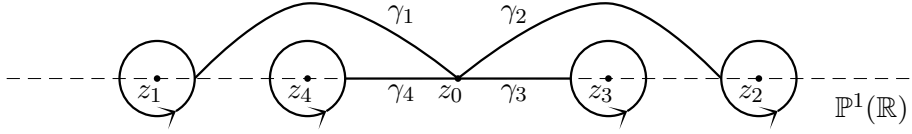
\begin{figure}[!ht]
\begin{center}
\begin{pspicture}(-6,-0.8)(6,1.5)
  \psline[linewidth=0.01,linestyle=dashed](-6,0)(6,0)
  \psdots(-4,0)(-2,0)(0,0)(2,0)(4,0) \rput(5.5,-0.4){$\PP^1(\RR)$}
  \rput(-0.1,-0.2){$z_0$} \rput(-4.1,-0.2){$z_1$}
  \rput(-2.1,-0.2){$z_4$} \rput(1.9,-0.2){$z_3$}
  \rput(3.9,-0.2){$z_2$}
  \pscurve[linewidth=0.03](0,0)(-2,1)(-3.5,0)
  \psarc[linewidth=0.03,arrowsize=8pt
  2,arrowlength=0.6,arrowinset=0.75]{->}(-4,0){0.5}{300}{300}
  \rput(-0.75,0.85){$\gamma_1$}
  \psarc[linewidth=0.03,arrowsize=8pt
  2,arrowlength=0.6,arrowinset=0.75]{->}(-2,0){0.5}{300}{300}
  \psline[linewidth=0.03](0,0)(-1.5,0) \rput(-0.75,-0.2){$\gamma_4$}
  \psarc[linewidth=0.03,arrowsize=8pt
  2,arrowlength=0.6,arrowinset=0.75]{->}(2,0){0.5}{300}{300}
  \psline[linewidth=0.03](0,0)(1.5,0) \rput(0.75,-0.2){$\gamma_3$}
  \pscurve[linewidth=0.03](0,0)(2,1)(3.5,0)
  \psarc[linewidth=0.03,arrowsize=8pt
  2,arrowlength=0.6,arrowinset=0.75]{->}(4,0){0.5}{300}{300}
  \rput(0.75,0.85){$\gamma_2$}
\end{pspicture}
\caption{\label{homotopic_basis}
  \textit{The standard homotopic basis}}
\end{center}
\end{figure}

\pbullet The complex conjugation also acts on the topological
fundamental group by left composition (thank you complex conjugation
for being continuous!). The fundamental group~$\pi_1(\PP^1_\CC
\setminus \{z_1,z_2,z_3,z_4\}, z_0)$ is simply denoted by~$\pi_1$. In
the rest of this paper, when we refer to the standard homotopic basis
of~$\pi_1$, we always mean the one drawn in the figure~\ref{homotopic_basis}.

\noindent So we have~$\gamma_1 \gamma_4 \gamma_3 \gamma_2 = 1$ and the complex conjugation acts
as follows:
\begin{equation} \label{gamma_gamma_bar}
\overline{\gamma_1} = \gamma_4^{-1} \gamma_1^{-1} \gamma_4
\qquad
\overline{\gamma_4} = \gamma_4^{-1}
\qquad
\overline{\gamma_3} = \gamma_3^{-1}
\qquad
\overline{\gamma_2} = \gamma_3 \gamma_2^{-1} \gamma_3^{-1}
\end{equation}

\pbullet Since complex conjugation is a continuous morphism of~$\CC$,
the monodromy of the cover~$\bar{p}$ can be deduced from the monodromy
of~$p$. More precisely, if we denote by~$\rho : \pi_1 \rightarrow
\Perm{F}$ and~$\bar{\rho} : \pi_1 \rightarrow \Perm{\bar{F}}$ the two
monodromy (anti)morphisms, then we have:
\begin{equation} \label{rho_rho_bar}
\forall \gamma \in \pi_1,
\qquad \bar{\rho}(\gamma) \circ c = c \circ \rho(\bar{\gamma})
\end{equation}

\pbullet From the Weil descent criterion (which boils down to the use
of Artin theory because the extension~$\CC/\RR$ is galois and finite)
we know that the cover~$p$ can be defined over~$\RR$ if and only if
there exists an isomorphism~$\Omega$ such that:
$$
\vcenter{ \xymatrix@=10pt{ \mathcal{C} \ar@{->}[rd]_{p}
    \ar@{->}[rr]^{\Omega} &
    & \bar{\mathcal{C}} \ar@{->}[ld]^{\bar{p}} \\
    & \PP^1_\CC }} \qquad \text{and} \qquad \bar{\Omega} \circ \Omega
= \Id
$$
(the last condition is a cocycle condition).

From these three points, we can deduce a completely combinatoric
criterion for the descent to~$\RR$ and for the existence of totally
real fibers:

\begin{theo} \label{criterion_tot_real}
  Let~$p : \mathcal{C} \rightarrow \PP^1_\CC$,~$z_0$, and~$F$ be as
  above.  Let~$(\sigma_1, \sigma_2, \sigma_3, \sigma_4) \in
  \Perm{n}^4$ be the branch cycle description of~$p$ ---~i.e.
  $\sigma_i = \rho(\gamma_i)$ where~$\rho : \pi_1 \rightarrow
  \mathcal{S}_{F}$ is the monodromy morphism~--- and by~$G =
  \ideng{\sigma_1, \sigma_2, \sigma_3, \sigma_4} \subset \Perm{F}$ the
  monodromy group of~$p$.
\begin{enumerate}
\item The cover~$p$ is defined over~$\RR$ if and only if there exists
  an involution~$\tau \in \Perm{F}$ such that:
  $$
  \sigma_4 \sigma_1^{-1} \sigma_4^{-1} = {}^\tau \sigma_1, \qquad
  \sigma_4^{-1} = {}^\tau \sigma_4, \qquad \sigma_3^{-1} = {}^\tau
  \sigma_3, \qquad \sigma_3^{-1} \sigma_2^{-1} \sigma_3 = {}^\tau
  \sigma_2
  $$
\item If that is the case and if moreover the cover~$p$ has no
  automorphism~--- i.e. if the centralizator of~$G$ in~$\Perm{F}$ is
  trivial~---, then the fiber~$F$ is totally real if and only if~$\tau
  =\Id$.
\end{enumerate}
\end{theo}

\noindent
{\slshape Proof}. --- Corresponding to the isomorphism~$\Omega$ of the
Weil descent criterion there is a bijection~$\omega : F \rightarrow
\bar{F}$; in term of~$\omega$ the conditions of the criterion are:
$$
\left[ \forall \gamma \in \pi_1, \quad \omega \circ \rho(\gamma) =
  \bar{\rho}(\gamma) \circ \omega \right] \quad \text{and} \quad
\left[ (c^{-1} \circ \omega)^2 = \Id \right]
$$
(the last condition is just the cocycle one). Let~$\tau = c^{-1}
\circ \omega$. This is an involution of~$\Perm{F}$ which
by~(\ref{rho_rho_bar}) and~(\ref{gamma_gamma_bar}) satisfies the
expected relations on the~$\sigma_i$ (be careful: the monodromy is an
anti-morphism) if and only if~$p$ is defined over~$\RR$.

Secondly, assuming that~$p$ is defined over~$\RR$, then the
isomorphism~$\Omega$ is an automorphism; so if moreover~$p$ has no
automorphism,~$\Omega$, and therefore~$\omega$, must be identity.
Furthermore, the conjugation~$c$ can now be viewed as a permutation
of~$F$.  In conclusion, the bijection~$\tau$ introduced in~$1$
satisfies~$\tau = c^{-1}$ and the fiber~$F$ is totally real if and
only if~$\tau = \Id$.\hfill $\square$

Having this result in mind, we can now go back to our family. Recall
that the covers of our family have no automorphism (since the center
of~$\Perm{n}$ is trivial). It is not difficult to verify that the
Nielsen class:
$$
a_{\frac{n}{2}-1} =
\begin{cases}
  \sigma_1 = (1, 2, \ldots, n-2) \\
  \sigma_2 = (n-2, n-1, n) \\
  \sigma_3 = \left[ \prod_{i = 1}^{\frac{n}{2} - 2} (i, n-i-2) \right] (n-2, n) (\frac{n}{2}-1)(n-1) \\
  \sigma_4 = \left[ \prod_{i = 1}^{\frac{n}{2} - 1} (i, n-i-1) \right] (n-1, n)\\
\end{cases}
$$
satisfies the previous theorem and is the only one in that case.

\begin{fact} \label{tot_real}
  Our family contains covers defined over~$\RR$ with an interval of
  non ramified points with totally real fibers.
\end{fact}

\subsection{Totally real~$\Perm{n}$-extensions. \normalfont ---}

At this point we know that our family contains some covers defined
over~$\QQ$ and some others with totally real fibers; we want to prove
that there are covers satisfying both properties.  We will have to
move one of the ramification point.

\reserve0/{\vbox{\hbox{
\begin{pspicture}(-1.5,-2)(1,2)
  \psdots(0,1.25)(0,-1.25)
  \pscircle(0,1.25){0.7} \rput(0.3,1.3){$h_4$} \rput(-1,1.5){$V$}
  \pscurve[linewidth=0.06](-0.6,1.61)(0,1.25)(0.4,0.675)
  \pscircle(0,-1.25){0.7} \rput(0.3,-1.1){$z_4$} \rput(-1,-1.5){$U$}
\pscurve[linewidth=0.06](-0.7,-1.1825)(0,-1.25)(0.6,-1.61) 
\psline{->}(0,0.5)(0,-0.5) \rput(0.5,0){$\varphi_{|V}$}
\rput(-1,0){$\genfrac{}{}{0pt}{}{\text{real}}{\text{points}}$}
\pscurve[linestyle=dotted]{->}(-1,0.3)(-0.8,1)(0,1.2)
\pscurve[linestyle=dotted]{->}(-1,-0.35)(-0.5,-0.5)(0,-1.2)
\end{pspicture}
}}} \endreserve10/ Suppose that three of the ramification points~$z_1,
z_2, z_3 \in \PP^1_\CC(\QQ)$ are fixed and let the fourth one~$z_4$
move on~$\PP^1(\RR)$ between~$z_1$ and~$z_3$. Thanks to the previous
section, for all such~$z_4$, there is a unique point~$h_4 \in
\varphi^{-1}(z_4)$ which represents a cover having with totally real
fibers above~$\left]z_2, z_3\right[$ (namely the cover with branch
cycle description~$a_{\frac{n}{2}-1}$ with respect to a standard
homotopic basis).  We choose~$U$ and~$V$ two neighborhoods of~$z_4$
and~$h_4$ respectively such that~$\varphi_{\mid V}$ becomes an
homeomorphism from~$V$ to~$U$. For every~$z \in U \cap \PP^1(\RR)$,
the covering corresponding to~$\varphi_{\mid V}^{-1}(z)$ satisfies the
preceding descent criteria; so it is defined over~$\RR$ and has a real
interval of specialization.  Thus we have~$\varphi_{\mid V}^{-1}(U
\cap \PP^1(\RR)) \subset V \cap \PP^1(\RR)$.  But rational points are
dense in~$\PP^1(\RR)$ so we can find rational points in~$V \cap
\PP^1(\RR)$; all the corresponding covers are defined over~$\QQ$ with
a complete segment of totally real fibers:

\begin{fact} \label{Q_plus_tot_real}
  Our family contains a rational a non empty interval of covers
  defined over~$\QQ$ each having an interval of totally real fibers.
\end{fact}

Let~$p : \PP^1_\QQ \rightarrow \PP^1_\QQ$ be such a cover.  Since an
interval is not a thin set, by Hilbert irreducible theorem, one can
find irreducible and totally real specialization.

\begin{prop} \label{existence_S_n}
  There exists totally real $\Perm{n}$-extensions of~$\QQ$ of
  degree~$n$.
\end{prop}

\subsection{Totally real~$\mathcal{A}_n$-extensions. \normalfont ---}

From the previous construction, we now want to deduce the same kind of
result for the group~$\mathcal{A}_n$. We need appropriate notations:
the cover~$p : \PP^1_\QQ \rightarrow \PP^1_\QQ$ will be denoted by~$p
: \mathcal{D} \rightarrow \mathcal{C}$ and its galois closure
by~$\mathcal{D}^{\rm gal} \rightarrow \mathcal{C}$. We also
consider~$\mathcal{C}' = \mathcal{D}^{\rm gal} / \mathcal{A}_n$
and~$\mathcal{D}'$ the irreducible component of the fiber
product~$\mathcal{D} \times_{\mathcal{C}} \mathcal{C}'$ which
satisfies the following diagram:
$$
\vcenter{ \xymatrix@R=3pt@C=8pt{
    \mathcal{D}^{\rm gal} \ar@{->}[dr] \\
    & \mathcal{D}' \\
    \mathcal{C}' \ar@{<-}[ur]_{n} \ar@{<-}[uu]^{\mathcal{A}_n} \\
    & \mathcal{D} \ar@{<-}[uu]_{2} \\
    \mathcal{C} \ar@{<-}[ur]^p_n \ar@{<-}[uu]^2 }} \qquad\qquad
\begin{array}{l}
\mathcal{C}' \stackrel{?}{\simeq}_\QQ \PP^1_{\QQ},
\quad \mathcal{D}' \stackrel{?}{\simeq}_\QQ \mathcal{D} \times_{\mathcal{C}} \mathcal{C}' \\
\\
\text{? means ``to be proved''}
\end{array}
$$
Since two of the inertia permutations are even and the two others
are odd, there are only two branched points in the cover~$\mathcal{C}'
\rightarrow \mathcal{C}$, i.e.~$z_1$ and~$z_3$ or~$z_4$ according
to~$n \equiv 0$ or~$2 \pmod{4}$. Therefore, by Riemann-Hurwitz
formula,~$\mathcal{C}'$ is a genus zero curve; since it has at least
one rational point (for example, one of the two ramified points), it
is also $\QQ$-isomorphic to~$\PP^1_\QQ$ (this kind of argument looks like the so-called ``double group trick'', see~\cite{topics}).

Now if we would have~$\QQ(\mathcal{C}') \subset \QQ(\mathcal{D})$,
then all the ramification indices of~$p$ above~$z_1$ would be
divisible by~$2$; it does not, so the fields~$\QQ(\mathcal{C}')$
and~$\QQ(\mathcal{D})$ are linearly disjoint over~$\QQ(\mathcal{C})$
and the curve~$\mathcal{D}'$ is nothing else than~$\mathcal{D}
\times_{\mathcal{C}} \mathcal{C}'$.

At this point, we have prove that the cover~$\mathcal{D}' \rightarrow
\mathcal{C}'$ has a Galois group equal to~$\mathcal{A}_n$. By
specialization, it gives degree~$n$ $\QQ$-extensions with Galois group
equal to~$\mathcal{A}_n$.  But we should also verify that there are
totally real fibers.  It suffices to show that the cover~$\mathcal{D}'
\rightarrow \mathcal{C}$ has totally real fibers. But,
since~$\mathcal{D}' \simeq \mathcal{D} \times_{\mathcal{C}}
\mathcal{C}'$, we can deduce the monodromy of the cover~$\mathcal{D}'
\rightarrow \mathcal{C}$ from these of the two covers~$\mathcal{D}
\rightarrow \mathcal{C}$ and~$\mathcal{C}' \rightarrow \mathcal{C}$.
Indeed, if~$(\sigma_1, \ldots, \sigma_4) \in (\Perm{n})^4$,
respectively\footnote{Two of the permutations~$\tau_i$ equal
  identity.}~$(\tau_1, \ldots, \tau_4) \in (\Perm{2})^4$, are the
branched cycle descriptions of the covers~$\mathcal{D} \rightarrow
\mathcal{C}$, respectively~$\mathcal{C}' \rightarrow \mathcal{C}$,
then:
$$
(\sigma_1 \times \tau_1, \ldots, \sigma_4 \times \tau_4) \in
(\Perm{\{1, \ldots, n\} \times \{1,2\}})^4
$$
is the branched cycle description of~$\mathcal{D}' \rightarrow
\mathcal{C}$. At last, the criterion~\ref{criterion_tot_real} shows
that the fiber above every point~$z_0$ between~$z_4$ and~$z_3$ is
totally real.

\begin{prop} \label{existence_A_n}
  There exists totally real $\mathcal{A}_n$-extensions of~$\QQ$ of
  degree~$n$.
\end{prop}

At the end of this paper, we give an explicit version of both
propositions~\ref{existence_S_n} and~\ref{existence_A_n}.

\section{Explicit computation}

In this section, we compute the Hurwitz space and the universal curve
for our family of covers. We fix once and for all three rational points~$z_1
< z_3 < z_2$. To ease notations we denote by~$\mathcal{H}$
the curve~${\mathcal H}'_{(z_1, z_2, z_3)}$.

\subsection{The universal curve and a choice of coordinates. \normalfont ---} \label{framework}

Because the covers in our family have no automorphism,
by~\cite{Fried_Volklein} or~\cite{Fried_Biggers}, there exists a
universal curve~$\mathcal S$ and a fibration:
$$
\xymatrix@R=18pt@C=0pt{ {\mathcal H} \ar[d]^{\varphi} &
  *[l]{\mathcal S}
  \ar[l]_{\Pi} \ar[d]^{\varPhi}\\
  \proC \smallsetminus \{z_1 , z_2 , z_3 \} & *[r]{F_2 (\proC
    \smallsetminus \{z_1 , z_2 , z_3 \} )} \ar[l]_{\qquad \pi} }
$$
where:
\begin{itemize}
\item $\mathcal S$ is a smooth quasi-projective surface over~$\QQ$,
\item $F_2 (\proC \smallsetminus \{z_1 , z_2 , z_3 \} )$ denotes the
  quasi-projective variety defined by the ordered pairs $(u,v) \in
  \proC \smallsetminus \{z_1 , z_2 , z_3 \}$ with~$u \neq v$,
\item $\pi$ is the morphism obtained by forgetting the second
  coordinate $v$,
\item the vertical arrows are finite and étale morphisms of varieties,
  all defined over~$\QQ$.
\end{itemize}
The morphism $\pi$ admits a well known projective completion that is
still denoted by $\pi$.  This is the canonical morphism from the
(projective) moduli space of curves of genus zero with $5$ marked
points $\mathcal M _{0,5}$ to the one with $4$ marked points $\mathcal
M _{0,4}$ (we refer to \cite{vanderput} for a highly comprehensive
study of the spaces $\mathcal M _{0,n}$ from an algebraic view point.
A lot of ideas contained in this paper are used here).  We
define~$\overline{\mathcal S}$ to be the normalization of~$\mathcal M
_{0,5}$ in~$\QQ(\mathcal{S})$.  We obtain this way a commutative
diagram between smooth projective varieties defined over $\QQ$ which
extends the previous one:
$$
\xymatrix{
  {\overline{\mathcal H}} \ar[d]^{\varphi} & \overline{\mathcal S} \ar[l]_{\Pi}\ar[d]_{\varPhi}\\
  \mathcal M_{0,4}& \mathcal M_{0,5} \ar[l]_{\pi} }
$$

In order to choose a system of coordinates on the universal curve, we
need to define another Hurwitz curve, a little bigger than ${\mathcal
  H}$, namely the moduli space of covers in our family with a marking
of an unramified point above~$z_3$.  This space denoted by~${\mathcal
  H}^{\bullet}$ is a degree~$2$ cover of~${\mathcal H}$ by the
forgetting map. We check that this cover is connected. There is a
universal curve~$\mathcal{S}^\bullet \rightarrow \mathcal{H}^\bullet$
obtained by base extension of~$\mathcal{S} \rightarrow \mathcal{H}$.
A normalization, as in the previous paragraph, provides the following
commutative diagram of projective varieties over $\QQ$:
$$
\xymatrix{
  {\overline{{\mathcal H}^{\bullet}}} \ar[d]^{\varphi} & \overline{{\mathcal S}^{\bullet}} \ar[l]_{\Pi^{\bullet}}\ar[d]_{\varPhi}\\
  \mathcal M_{0,4}& \mathcal M_{0,5} \ar[l]_{\pi} }
$$

Now, we want to choose adapted coordinates on those spaces. To begin
with, let's recall that there exist four sections $s_1$, $s_2$, $s_3$
and $s_4$ of $\pi$, corresponding to the four marked points (see
\cite{vanderput} section $3$). We define $x$ and $y$ to be the
coordinates on $\mathcal M_{0,5}$ such that $x=y=\infty$ on~$s_1$,
$x=1$ on~$s_2$, $x=y=0$ on~$s_4$ and $y=1$ on~$s_3$. Then we
set~$\lambda \affecte \frac{x}{y} \in \CC ({\mathcal M}_{0,4})$.

We also need coordinates on the universal curve~$\overline{{\mathcal
    S}^{\bullet}}$.  The fibration $\overline{{\mathcal S}^{\bullet}}
\rightarrow {\overline{{\mathcal H}^{\bullet}}}$ admits four sections
corresponding to the four points~$A,B,C,D$. This amounts to saying that
these four points are defined over~$\CC(\overline{{\mathcal
    H}^{\bullet}})$. This is clear for~$A$ and~$B$ because they are
isolated.  This is also true for~$C$ and~$D$ by definition
of~$\overline{\mathcal{S}^\bullet}$.  We define the function~$X$ to be
the unique coordinate taking values~$\infty, 1, 0$ at~$A,B,C$
respectively. We define~$Y$ to be the unique coordinate taking
values~$\infty, 0, 1$ at~$A,C,D$ respectively. We set~$\mu \affecte
\frac{X}{Y} \in \CC(\overline{{\mathcal H}^{\bullet}})$.  This
situation is summarized in table~\ref{pointing}.

\begin{table}[!ht]
{\scriptsize
  $$
  \xymatrix@R=10pt@C=0pt{
    \infty  &         & \infty &     &     & 1     &     &         & 0    &       & 0   & 1     &                   & &             &    \\
    & (A,n-2) \ar@{->}[ul]^X \ar@{->}[ur]_Y  &        & & (\bullet, 1) & (B,3) \ar@{->}[u]^X &     & (\bullet, 1)^{n-3} &      & (C,1) \ar@{->}[ul]^X \ar@{->}[ur]_Y &     & (D,1) \ar@{->}[u]_Y & (\bullet, 2)^{\frac{n}{2}-1} & & (\bullet, 2)^{\frac{n}{2}} & \\
    \\
    &         &        & z_1 \ar@{-}[uull] \ar@{=}[uur] &        &       & z_2 \ar@{-}[uul] \ar@{=}[uur] &                   &      &       &         & z_3 \ar@{-}[uull]\ar@{-}[uu]\ar@{=}[uur]  &    &              &  z_4  \ar@{=}[uu] &  \\
    & & \infty \ar@{<-}[ur]^x & & \infty \ar@{<-}[ul]_y & 1
    \ar@{<-}[ur]^x & & 1/\lambda \ar@{<-}[ul]_y & & & \lambda
    \ar@{<-}[ur]^x & & 1 \ar@{<-}[ul]_y & 0 \ar@{<-}[ur]^x & & 0
    \ar@{<-}[ul]_y }
  $$
  }
\caption{\label{pointing} Pointing the covers of ${\mathcal H}$ and choice of coordinates}
\end{table}

In order to compute an algebraic model for the cover $\varPhi$, we
first exhibit an explicit model for a degenerate cover. In other words
we compute the residual morphism induced by $\varPhi$ on the fiber
over a boundary point of $\mathcal M _{0,4}$ (that is corresponding to
$z_1$, $z_2$ or $z_3$).  Next, we can rebuilt the entire family by
formal deformation of this residual morphism viewed as a morphism of
curves on a formal power series ring.

\subsection{Degenerate covers and their computation. \normalfont ---} \label{s_revts_degeneres}
Let $b$ be a boundary point of $\mathcal M _{0,4}$. We fix a point,
e.g. $b = z_3$. Since the key point of this section is the explicit
algebraic structure of the compactification of moduli spaces of curves
of genus zero, we will assume that the reader has some familiarities
with these notions. We just recall (see~\cite{vanderput}, section $3$)
that the fiber~$\mathcal{C}_b = \pi^{-1}(b)$ is a stable $4$-pointed
tree of projective lines made of two irreducible components. Let us
denote by ${\mathcal C}_{b,1}$ the irreducible component
of~$\mathcal{C}_b$ containing the two closed points $z_1$ and $z_2$,
and ${\mathcal C}_{b,2}$ the other one.  For each point~$h \in
\overline{\mathcal{H}^\bullet}$ such that $\varphi (h) = b$,
set~$\mathcal{D}_h = {(\Pi^\bullet)}^{-1}(h)$. Then the
restriction~$\varPhi_h$ of~$\varPhi$ on~$\mathcal{D}_h$ is a cover of
nodal curves (see for example figure~\ref{revt_degenere}). We denote
by~$\varPhi_{h,i}$ the restriction~$\varPhi_h$ to~$\mathcal{D}_{h,i} =
\varPhi_h^{-1}(\mathcal{C}_i) \cap \mathcal{D}_h$.

For $i=1,2$ the morphism $\varPhi_{h,i}$ is finite and the
ramification locus is contained in the union of the two marked points
and the singular point. Now the monodromy can classically be deduced
(see \cite{boundary} for example) from the one of the non-degenerate
covers in a small neighborhood of $h$. More precisely, let~$V$ be a
small enough neighborhood of $h$ (for the complex topology).  If the
branch cycle description of covers parameterized by~$V \setminus
\{h\}$ is given by a $4$-tuple $[ \sigma_1 , \sigma_2 , \sigma_3 ,
\sigma_4]$ then the branch cycle description of $\varPhi_{h,1}$
(respectively~$\varPhi_{h,2}$) is $[ \sigma_1\sigma_2 , \sigma_3 ,
\sigma_4]$ (respectively~$[ \sigma_1 , \sigma_2 , \sigma_3\sigma_4]$).

Let us now concentrate on the computation of an algebraic model for
the covers $\varPhi_{h,1}$ and $\varPhi_{h,2}$. Recall, from the
beginning of this section, that we have fixed three rational
points~$z_1 < z_3 < z_2$.  For any~$z_4 \in \left]z_1,z_3\right[$
let~$p_{z_4}$ be the cover of~$\PP^1 \setminus \{z_1,z_4,z_3,z_2\}$
with monodromy~$a_{\frac{n}{2}-1}$ in the homotopic basis represented
in figure~\ref{homotopic_basis}.  Letting~$z_4$ tends to~$z_3$ we
define a point~$h$ in the boundary
of~$\overline{\mathcal{H}^\bullet}$. For this~$h$, the ramification
data for $\varPhi_{h,2}$ is: {\scriptsize
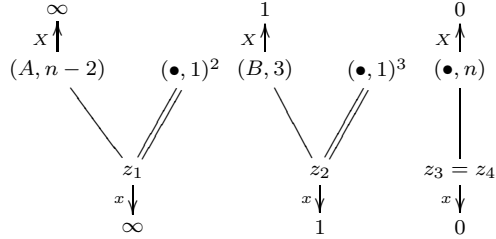
  $$
  \xymatrix@R=10pt@C=0pt{
    \infty    \ar@{<-}[d]_X & & & 1 \ar@{<-}[d]_X & & & 0 \ar@{<-}[d]_X \\
    (A, n-2) & & (\bullet, 1)^2  & (B,3) & & (\bullet, 1)^3 & (\bullet, n) \\
    \\
    & z_1 \ar@{-}[uul] \ar@{=}[uur] & & & z_2 \ar@{-}[uul] \ar@{=}[uur] & & z_3 = z_4 \ar@{-}[uu] \\
    & \infty \ar@{<-}[u]^x & & & 1 \ar@{<-}[u]^x & & 0 \ar@{<-}[u]^x }
  $$} This cover is a Padé approximant and for the couple of
coordinates chosen, we have:
\begin{equation}\label{Phi_h_2}
x = \varPhi_{h,2} (X) = \frac{X^n}{\frac{n(n-1)}{2} \left(X^2 - 2
    \frac{n-2}{n-1} X + \frac{n-2}{n}\right)}
\end{equation}
Similarly for $\varPhi_{h,1}$, we find:
\begin{equation}\label{Phi_h_1}
y = \varPhi_{h,1} (Y) = \frac{1}{2} T_n(2Y-1) + \frac{1}{2}
\end{equation}
where $T_n$ denotes the Tchebycheff polynomial of order $n$.

Therefore a very simple model for this degenerate cover is known for
every even~$n$. We stress this decisive fact for the
relevance of our approach.

\subsection{Effective deformation and formal patching. \normalfont ---}

In order to built the entire family, we now formally deform the
previous degenerate cover.  The deformation technique for covers
appeared in~\cite{Fulton_Hur} and were then developed by D.~Harbater
for the study of the inverse Galois problem over complete local fields
(see e.g.~\cite{Harbater_deformation, Harbater_GalCov}).
In~\cite{Wewers_Def}, S.~Wewers gives a presentation of the technique
of deformation very well adapted to our purpose.

\subsubsection{Computation of a formal model using effective deformation. \normalfont ---}

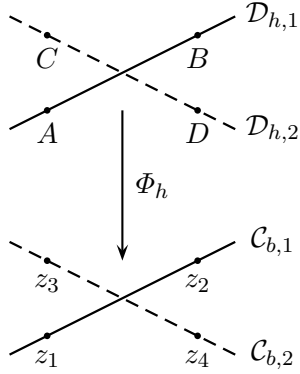
\begin{figure}[!ht]
  \centerline{
\begin{pspicture}(-1.6,-1)(2.2,4)
  \psline[linewidth=0.03](-1.5,2.25)(1.5,3.75)
  \psline[linewidth=0.03,linestyle=dashed](-1.5,3.75)(1.5,2.25)
  \psdots(-1,2.5)(-1,3.5)(1,2.5)(1,3.5)
  \rput(-1,2.2){$A$}\rput(-1,3.2){$C$}\rput(1,2.2){$D$}\rput(1,3.2){$B$}
  \rput(2,3.75){$\mathcal{D}_{h,1}$}\rput(2,2.25){$\mathcal{D}_{h,2}$}
  \psline[linewidth=0.03](-1.5,-0.75)(1.5,0.75)
  \psline[linewidth=0.03,linestyle=dashed](-1.5,0.75)(1.5,-0.75)
  \psdots(-1,-0.5)(-1,0.5)(1,-0.5)(1,0.5)
  \rput(-1,-0.8){$z_1$}\rput(-1,0.2){$z_3$}\rput(1,-0.8){$z_4$}\rput(1,0.2){$z_2$}
  \rput(2,0.75){$\mathcal{C}_{b,1}$}\rput(2,-0.75){$\mathcal{C}_{b,2}$}
  \psline{->}(0,2.5)(0,0.5)\rput(0.4,1.5){$\varPhi_h$}
\end{pspicture}
}
\caption{\label{revt_degenere}
  \textit{The degenerate cover~$\varPhi_h$}}
\end{figure}

Let us denote by~$R$ the complete local ring
of~$\overline{\mathcal{H}^\bullet}$ at the point~$h$,
namely~$\CC[[\mu]]$ because~$\mu$ is a local parameter at~$h$. By base
extension~$\mathcal{S}_R = \overline{\mathcal{S}^\bullet}
\times_{\overline{\mathcal{H}^\bullet}} \Spec(R)$ and~$\mathcal{C}_R =
\mathcal{M}_{0,5} \times_{\mathcal{M}_{0,4}} \Spec(R)$ are {\em
  projective nodal curves} whose special fibers are nothing else
than~$\mathcal{D}_h$ and~$\mathcal{C}_b$ respectively. The
cover~$\varPhi$ induces a {\em tame admissible cover}
from~$\mathcal{S}_R$ to~$\mathcal{C}_R$ which is a deformation of the
previous degenerate cover~$\varPhi_h$ represented in
figure~\ref{revt_degenere}.

We describe the {\em deformation datum} associated to our deformation
as it is explained in the paragraph of~\cite{Wewers_Def} entitled {\em
  The main result} (pages~240-241). Using the preceding and putting~$z
= 1/y$, the curve~$\mathcal{C}_R$ has affine equation~$xz = \lambda$.
So there is only one ordinary double point~$\delta$ i.e.~$(0:0:1)$.
Moreover, the complete local ring~$\OO_{\mathcal{C}_R, \delta} \simeq
R[[x,z]] / \ideng{xz - \lambda}$.  We also have a {\em mark}
on~$\mathcal{C}_R$, namely the horizontal divisor defined by the four
generic branched points
$$
(0:1:0), \quad (1:\lambda:1), \quad (1:0:0), \quad (\lambda:1:1)
$$
and this divisor is étale over~$\Spec(R)$. The
curve~$\mathcal{S}_R$ has also a unique singular point~$\Delta$ whose
local ring is isomorphic to~$R[[X,Z]] / \ideng{XZ - \mu}$ where~$Z =
1/Y$.  Moreover, due to the ramification of the degenerate cover, we
know that~$X^n \sim x$, $Z^n \sim z$ and~$\mu^n \sim \lambda$
(where~$\sim$ means equal up to a unit factor). Therefore the {\em
  deformation datum} consists only in~$\mu = XZ \in R$.

In concrete terms, effective patching permits us to compute a model of
our family over a Puiseux series field~$\CC((\mu))$ where this field
is nothing else than the completion
of~$\CC(\overline{\mathcal{H}^\bullet})$ at the point~$h$.  We will
call this model an formal one.

In view of our choice of coordinates, the model we are looking for has
the following form:
\begin{align*}
  &S(X) = \frac{(X^{\frac{n}{2}} + \alpha_{\frac{n}{2} - 1}
    X^{\frac{n}{2} - 1} + \cdots + \alpha_0)^2}
  {\gamma (X^2 + \beta_1 X + \beta_0)} = \frac{S_0(X)}{S_\infty(X)} \\
  &S(X) - 1 = \frac{(X - 1)^3(X^{n-3} + \delta_{n-4} X^{n-4} + \cdots
    + \delta_0)}
  {\gamma (X^2 + \beta_1 X + \beta_0)} = \frac{S_1(X)}{S_\infty(X)} \\
  &S(X) - \lambda = \frac{X(X - \varepsilon_0)(X^{\frac{n}{2} - 1} +
    \eta_{\frac{n}{2} - 2} X^{\frac{n}{2} - 2} + \cdots + \eta_0)^2}
  {\gamma (X^2 + \beta_1 X + \beta_0)} =
  \frac{S_\lambda(X)}{S_\infty(X)}
\end{align*}
where the~$2n$ coefficients~$\alpha_i, \beta_j, \delta_k,
\eta_l$,~$\varepsilon_0$ and~$\gamma$ have to be found
in~$\CC(\overline{\mathcal{H}^\bullet})$. From the three equalities
above, one can deduce that~$S_0(X) - S_1(X) - S_\infty(X) = 0$ and
that~$S_0(X) - S_\lambda(X) - \lambda S_\infty(X) = 0$. This gives us
a system of~$2n$ polynomials in~$2n$ variables and with coefficients
in~$\CC[\mu]$ (because~$\lambda$ is a polynomial in~$\mu$) which
should be satisfied by our $2n$ coefficients.

The knowledge of the degenerate cover gives us first order $\mu$-adic
development for all the coefficients. Then computing higher orders is
just a careful application of the Newton-Hensel lemma.

\subsubsection{Computation of an algebraic model from the formal one. \normalfont ---} First of all, we change the coordinate~$X$ by an homography which fixes~$1$ and~$\infty$ so as to cancel the
coefficient~$\alpha_{\frac{n}{2}-1}$ in the polynomial~$S_0(X)$.
$$
X \longleftarrow \left(1 + \frac{\alpha_{\frac{n}{2} -
      1}}{\frac{n}{2}}\right) X - \frac{\alpha_{\frac{n}{2} -
    1}}{\frac{n}{2}}
$$
This natural normalization turns all the coefficients in our model
to be defined over~$\CC(\mathcal{H})$.  Hopefully, we have noticed the
all the new coefficients are now power series in~$\mu^2$; so we are
well back on~$\mathcal{H}$!

Then the last part of the computation consists in deriving an
algebraic model over~$\CC(\mathcal{H})$ from the preceding formal one
defined over~$\CC((\mu^2))$ (which, we recall; is the completion
of~$\CC(\mathcal{H})$ with respect to a point of~$\mathcal{H}$).
Theoretically speaking, this steps is based on the Artin's
algebraization theorem. From a computational point of view, using the
known $\mu$-adic approximations of the coefficients, we should:

\smallbreak

\pbullet first deduce a generator of~$\CC(\mathcal{H})$;

\pbullet and then express all the coefficients as rational fractions
in this generator.

\smallbreak

Even if the first step can be done systematically (as is explained
in~\cite{JMC_Tools}), we just guess a generator~$T$ among the
coefficients and compute all the other coefficients in function
of~$T$. If we know the $\mu$-adic approximations of~$T$ and of every
coefficients~$C$ with enough accuracy, finding such expression is just
a matter of linear algebra; indeed, for increasing values of the
degree~$d$, we have to solve in~$\alpha_i, \beta_j \in \CC$ an
equation like:
$$
\alpha_d T^d + \cdots + \alpha_0 + C (\beta_d T^d + \cdots +
\beta_0) = 0
$$
which expanded in~$\mu$ gives rise to a linear system in~$\alpha_i,
\beta_j$.  We stress that \cite{JMC_Tools} gives an upper bound for
the degree~$d$

\begin{rk}
  The computation of the last two steps could involve computations
  with complex numbers, and so numerical approximations, because
  nothing tells us that the model we have chosen is defined
  over~$\QQ$. But, luckily, it is!
\end{rk}


\section{Numerical results}

The two covers we are looking for are given by:
$$
\begin{array}{rccc}
\varphi : & \mathcal{H} & \longrightarrow & \mathcal{M}_{0,4} \\
       &     T       & \longmapsto     & H_n(T)
\end{array}
\qquad\qquad
\begin{array}{rccc}
\varPhi : & \mathcal{S} & \longrightarrow & \mathcal{M}_{0,5} \\
       &    X        & \longmapsto     &   S_n(T, X)
\end{array}
$$
For~$n = 6$, we obtain: {\small
\begin{align*}
  &S_6(T, X) = \frac{\left( X^3 + \frac{75T + 120}{16}X + \frac{625T^3
        + 3600T^2 + 6720T + 4096}{96T + 256} \right)^2}
  {\frac{3(25T+56)^3}{2^8(3T+8)} \left(X^2 + TX + \frac{25T^3 + 120T^2
        + 192T + 128}{36T + 96}\right)
    } \\
  &S_6(T, X) - 1 = \frac{(X-1)^3 \left( X^3 + 3X^2 + \frac{75T +
        168}{8}X + \frac{625T^3 + 4950T^2 + 12960T + 11136}{48T + 128}
    \right)} {\frac{3(25T+56)^3}{2^8(3T+8)}
    \left(X^2 + TX + \frac{25T^3 + 120T^2 + 192T + 128}{36T + 96}\right)} \\
  &S_6(T, X) - H_6(T) = \\
  &\frac{ \left( X^2 - \frac{5T + 8}{2}X + \frac{125T^3 + 1050T^2 +
        2720T + 2176}{48T + 128} \right) \left( X^2 + \frac{5T +
        8}{4}X + \frac{25T^3 + 180T^2 + 424T + 320}{24T + 64}
    \right)^2 } {\frac{3(25T+56)^3}{2^8(3T+8)} \left(X^2 + TX +
      \frac{25T^3 + 120T^2 + 192T + 128}{36T + 96}\right)}
\end{align*}
} and:
$$
H_6(T) = \frac{(T+8)(T+\frac{13}{5})^2(T+\frac{8}{5})^3} {-3 \times
  5(T+\frac{8}{3})(T+\frac{56}{25})^3} \quad H_6(T) - 1 =
\frac{(T+2)(T+\frac{16}{5})^5} {-3 \times
  5(T+\frac{8}{3})(T+\frac{56}{25})^3}
$$

We manage similar computation in {\scshape Magma} for all $n$ up
to~$20$ in less than 20 minutes on an AMD 700Mhz.

In order to obtain totally real $\QQ$-extension of degree~$n$ with
Galois group~$\Perm{n}$ and~$\mathcal{A}_n$, we have to specialize
twice.  The method points out a special value~$t_h \in \QQ$ of
the parameter corresponding to the cover~$h$ we have deformed.
Let us choose a close enough rational number~$t_0 < t_h$ and set~$t
= t_0$.  We get a regular $\Perm{n}$-extension:
$$
\vcenter{ \xymatrix@R=3pt@C=8pt{
    L = \QQ(X)^{\rm gal} \ar@{-}[dr] \\
    & \QQ(X) \\
    \QQ(x) \ar@{-}[ur]_{n} \ar@{-}[uu]^{\mathcal{S}_n} }} \qquad
\qquad \text{where:} \quad x = S_n(t_0, X) \in \QQ(X)
$$
with four ramified points~$0, H_n(t_0), 1$ and~$\infty$ in this
order.  Moreover, due to the criterion of
section~\ref{s_totally_real_fibers}, all the points in~$\left]0,
  H_n(t_0)\right[$ have totally real fibers. By Hilbert
irreducibility theorem, most rationals in this interval provide totally real
$\mathcal S_n$-extensions of $\QQ$. The case of $\mathcal A_n$ follows classically.

\end{document}